\documentclass[12pt]{amsart}
\theoremstyle{plain}
\newtheorem{Thm}{Theorem}
\newtheorem{Ass}[Thm]{Assertion}

\newtheorem{Def}[Thm]{Definition}

\newtheorem{Main}{Main Theorem}

\begin{document}

\title[Ricci expanding solitons]
{Remarks on non-compact complete Ricci expanding solitons }

\author{Li MA and Dezhong CHEN}

\address{Department of mathematical sciences \\
Tsinghua university \\
Beijing 100084 \\
China}

\email{lma@math.tsinghua.edu.cn} \dedicatory{}
\date{July. 20th, 2005}

\thanks{$^*$ This work is supported in part by
the Key 973 project of Ministry of Science and Technology of China.
}

\begin{abstract}
In this paper, we study  gradient Ricci expanding solitons $(X,g)$
satisfying
$$
Rc=cg+D^2f,
$$
where $Rc$ is the Ricci curvature, $c<0$ is a constant, and $D^2f$
is the Hessian of the potential function $f$ on $X$.
 We show that for a gradient expanding soliton $(X,g)$ with
non-negative Ricci curvature, the scalar curvature $R$ has at least
one maximum point on $X$, which is the only minimum point of the
potential function $f$. Furthermore, $R>0$ on $X$ unless $(X,g)$ is
Ricci flat. We also show that there is exponentially decay for
scalar curvature for $\epsilon$-pinched complete non-compact
 expanding solitons.
\end{abstract}

 \maketitle

\section{Introduction}

In this paper, we continue our study on Ricci solitons \cite{M04},
which are generated by one parameter family of diffeomorphisms and
are special solutions to Ricci flow introduced by R.Hamilton in 1982
\cite{H095}. We assume in this paper that $(X,g)$ is a gradient
expanding soliton. Let recall the definition of expanding
 soliton.

\begin{Def}
 We call a Riemannian manifold $(X,g)$ an expanding soliton if
there is a smooth
  solution $f$ on a Riemannian manifold $(X,g)$ such that for some constant $c<0$, it holds the equation
  \begin{equation}
Rc=cg+D^2f, \label{Rc}
  \end{equation}
  on $X$,
  where $D^2f$ is the Hessian matrix of the function $f$ and $Rc$
  is the Ricci tensor of the metric $g$. We call the function $f$
  the potential function for the soliton $(X,g)$.
 If $c>0$ in (\ref{Rc}),  $(X,g)$ is called a shrinking soliton; if $c=0$, $(X,g)$
  is called a steady
  soliton.
\end{Def}

In the study of Ricci flow, we often meet the following definition.

\begin{Def}
 The Ricci curvature of a Riemannian manifold $(X,g)$ is called
$\epsilon$-pinched if there is some $\epsilon>0$ such that the
scalar curvature $R>0$ on $X$ and
$$
Rc\geq \epsilon Rg
$$
on $X$.
\end{Def}

  Throughout this paper, we shall assume that the Riemannian manifold
 $(X,g)$ is a complete non-compact Riemannian manifold of dimension $n\geq 3$.
We denote by $R$ the scalar curvature of the metric $g$.

Our main result is the following
\begin{Main}
Assume that the Ricci curvature of the gradient expanding soliton
$(X,g)$ is non-negative. Then the scalar curvature $R$ has at least
one maximum point on $X$, which is the only minimum point of the
potential function $f$. Furthermore, $R>0$ on $X$ unless $(X,g)$ is
Ricci flat.
\end{Main}

The proof of this Theorem will be proved in section 3.

 In section four, we will prove the following result

 \begin{Thm} Assume that $(X,g)$ is a gradient expanding soliton
with its Ricci curvature being $\epsilon$-pinched. Then its scalar
curvature has the decay
$$
R(s)\leq R(o)e^{Cs-Cs^2}.
$$
as the distance function $s$ from a fixed point going to infinity,
i.e., $s=d(x,o)\to+\infty$.
 \end{Thm}

 We remark that a similar
but weaker decay result is announced by L.Ni in Proposition 3.1 in
\cite{N05}. We know the result for a while and a reason for the
delay of this present is that we try to prove non-existence of this
kind of expanding solitons. However, we have not been succeed yet.

Throughout $C$ will denote various uniform constants in different
places.

\section{preliminary}
We recall first some basic properties about Ricci solitons
\cite{H95}.

Taking the trace of both sides of
 (\ref{Rc}), we have
 \begin{equation}
R=nc+\Delta f. \label{Delt}
 \end{equation}

Take a point $x\in X$. In local normal coordinates $(x^i)$ of the
Riemannian manifold $(X,g)$ at a point $x$, we write the metric $g$
as $(g_{ij})$. The corresponding Riemannian curvature tensor and
Ricci tensor are denoted by $Rm=(R_{ijkl})$ and $Rc=(R_{ij})$
respectively. Hence,
$$
R_{ij}=g^{kl}R_{ikjl}
$$
and
$$
R=g^{ij}R_{ij}.
$$
We write the covariant derivative of a smooth function $f$ by
$Df=(f_i)$, and denote  the Hessian matrix of the function $f$ by
$D^2f=(f_{ij})$, where  $D$ the covariant derivative of $g$ on $X$.
 The
higher order covariant derivatives are denoted by $f_{ijk}$, etc.
Similarly, we use the $T_{ij,k}$ to denote the covariant derivative
of the tensor $(T_{ij})$. We write $T_j^i=g^{ik}T_{jk}$. Then the
Ricci soliton equation is
$$
R_{ij}=f_{ij}+cg_{ij}.
$$
Taking covariant derivative, we get
$$
f_{ijk}=R_{ij,k}.
$$
So we have
$$
f_{ijk}-f_{ikj}=R_{ij,k}-R_{ik,j}.
$$
By the Ricci formula we have that
$$
f_{ijk}-f_{ikj}=R_{ijk}^lf_l.
$$
Hence we obtain that
$$
R_{ij,k}-R_{ik,j}=R_{ijk}^lf_l.
$$
Recall that the contracted Bianchi identity is
$$
R_{ij,j}=\frac{1}{2}R_i.
$$
Upon taking the trace of the previous equation we get that
$$
\frac{1}{2}R_i+R_{i}^kf_k=0,
$$
i.e.,
\begin{equation}
R_k=-2R_{k}^jf_j. \label{RR}
\end{equation}
Then at $x$,
$$
D_k(|Df|^2+R+2cf)=2f_j(f_{jk}-R_{jk}+2cg_{jk})=0.
$$
So,
\begin{equation}
|Df|^2+R+2cf=M, \label{MM}
\end{equation}
where $M$ is a constant.

In the remaining part of this section, we assume that $0\leq Rc\leq
C$ on the expanding soliton $(X,g)$ for some constant $C>0$. Then we
have $|D^2f|\leq C$ on $X$. Assume $f\geq 0$ and that $o$ is a
critical point of the potential function $f$. Then using the
Taylor's expansion, we have
$$
f(x)\leq C d^2(x,o). \label{B3}
$$

We now study the behavior of the potential function along a
minimizing geodesic curve on the expanding soliton. A similar work
has been done by G.Perelman \cite{P03} (see also \cite{KL}) where he
tries to to give some uniform bounds on potential function $f$ on a
shrinking soliton. Fix a point $o\in X$. Take any minimizing
geodesic curve $\gamma(s)$ connecting $x$ and the fixed point $p$,
where $s$ is the arc-length parameter. Write by $r=d(x,o)$ and
$X=\gamma'(s)$. Assume that $r>2$. Let $\{Y_i\}$ ($i=1,...n-1$) be
an orthonormal parallel vector fields along $\gamma$. Let $Y$ be an
orthogonal vector field along the curve $\gamma$ vanishing at end
points. Then the second variational formula  \cite{SY94} (see also
\cite{Aub82}) tells us that
$$
\int_0^r(|Y|^2-<R(X,Y)Y,X))ds\geq 0.
$$
Take $Y$ to be $sY_i$ on $[0,1]$, $=Y_i$ on $[1,r-r_0]$ where
$1<r_0<r$, and $\frac{r-s}{r_0}Y_i$. Adding over $i$ gives that
$$
\int_0^{r}Rc(X,X)\leq
C_0(r_0)+\frac{n-1}{r_0}-\int_{r-r_0}^r(\frac{r-r_0}{r_0})^2Rc(X,X)ds,
$$
which implies that for some constant $C>0$,
\begin{equation}
\int_0^{r}Rc(X,X)\leq C. \label{B1}
\end{equation}
Note that
$$
(\int_0^{r}Rc(X,Y_1)ds)^2\leq r\int_0^{r}|Rc(X,Y_1)|^2ds\leq
s\sum_i\int_0^{r}|Rc(X,Y_i)|^2ds.
$$
Thinking of $Rc$ as self-adjoint linear operator on $TX$ and taking
a point-wise orthonormal frame $\{e_j\}$ as eigenvectors of
$Rc=(\bigoplus\lambda_j)$ , we have that
$$
R=\sum_j\lambda_j
$$
and for $X=\sum_jX_je_j$,

$$
\sum_i|Rc(X,Y_i)|^2=<X,Rc^2X>=\sum_j\lambda_jX_j^2\leq RRc(X,X).
$$
Then,
$$
(\int_0^{r}Rc(X,Y_1)ds)^2\leq Cs \int_0^rRc(X,X)\leq C^2s.
$$
Hence, for any unit vector field $Y$ along $\gamma$, orthogonal to
$X$, we have
$$
\int_0^{r}Rc(X,Y)ds)\leq C(\sqrt{s}+1).\label{B2}
$$
Using (\ref{Rc}) we have
$$
\frac{d^2f(\gamma(s))}{ds^2}=Rc(X,X)-c\geq -c,
$$
and
$$
\frac{d(Yf)(\gamma(s))}{ds}=Rc(X,Y).
$$
Then we have
$$
\frac{df(\gamma(s))}{ds}\geq \frac{df(\gamma(s))}{ds}(0) -cs\geq
-cs+C
$$
and for $s>2$,
\begin{equation}
|(Yf)(\gamma(s))|\leq |(Yf)(\gamma(0))|+\int_0^s|Rc(X,Y)|ds\leq
C\sqrt{s}. \label{B5}
\end{equation}
Therefore, we can conclude that at large distance from $o$ the
potential function $f$ has its gradient making small angle with the
gradient of the distance function from $o$.

\section{Proof of Main Theorem}

Assume that $Rc\geq 0$ on $X$,  and we also assume that for some
constant $c<0$ we have $Rc-cg>0$ on $X$. By (\ref{Rc}) we know that
$$
D^2f=Rc-cg\geq -cg>0, \;\;\;on\;\;\;X.
$$
 Then the potential function $f$ is locally strictly convex. Since
 $(X,g)$ is a complete non-compact Riemannian manifold, we have that
 $f$ has at most one critical point,i.e., the point where $\nabla
 f=0$.
 Using $D^2f>0$, we know that if $p\in X$ is the critical point of
 $f$, then it is a non-degenerate minimum point of $f$.

Note that along any minimizing geodesic curve $\gamma(s)$ connecting
$x$ and the fixed point $p$, where $s$ is the arc-length parameter,
we have
\begin{align}
<\nabla f,
\gamma'(s)>|_0^s&=\int_0^sf_{ij}\frac{d\gamma^i}{ds}\frac{d\gamma^j}{ds}ds
\label{In}
\\
&=\int_0^s(R_{ij}-cg_{ij})\frac{d\gamma^i}{ds}\frac{d\gamma^j}{ds}ds\nonumber\\
&=-cs+\int_0^sR_{ij}\frac{d\gamma^i}{ds}\frac{d\gamma^j}{ds}ds\nonumber\\
&\geq -cs>0\nonumber
\end{align}
This implies that $f(\gamma(s))$ is growing at infinity at least the
quadratic rate $-c$ of the distance function. Then $f$ has at least
a minimum point in $X$.

Assume that $o$ is the only critical point of $f$.
 Then by adding a constant, we can
 assume that and $f(o)=0$ and $f>0$ on $X-\{o\}$. Using (\ref{MM}),
 we know that
 $$
M=|Df|^2(o)+R(o)+2cf(o)=R(o).
 $$
Using (\ref{RR}) we know that $o$ is also the critical point of $R$.

Let $x\in X-\{o\}$. Taking a minimizing geodesic curve $\gamma(s)$
connecting $x$ and the fixed point $o$, where $s$ is the arc-length
parameter, we again have by using $(\ref{In})$
$$
<\nabla f, \gamma'(s)>>-cs>0.
$$
This implies that the integral curves of $\nabla f$ in $X-\{o\}$
emanating from the point $o$ to infinity. Take a integral curve
$\sigma(t)$ $\nabla f$ in $X-\{o\}$. Then by (\ref{RR}) we have
\begin{equation}
\frac{d}{dt}R(\sigma(t)=R_if_i=-2Rc(\nabla f,\nabla f)\leq
0.\label{De}
\end{equation}
Hence $R(x)\leq R(o)$ for all $x\in X\{o\}$. So, $o$ is a maximum
point of $R$.

By this we conclude that
\begin{Ass} Assume that the Ricci curvature of the gradient expanding soliton $(X,g)$ is
non-negative positive. Then the scalar curvature $R$ has at least
one maximum point of $R$, which is the only critical point of the
potential function $f$.
\end{Ass}

If $R(o)=0$, then $R=0$ on $X$. Hence $Rc=0$ on $X$, that is to say
that $(X,g)$ is Ricci flat. So we have $R(o)>0$. By the local strong
maximum principle, we must have $R>0$ on the whole space $X$.

This finishes the \emph{proof of Main Theorem}.

In the remaining part of this section, we consider the behavior of
$f$ at infinity.
 Since
$$
|Df|(x)^2+2cf(x)=R(o)-R(x)\geq 0,
$$
we get that
$$
|Df|^2\geq -2cf=2|c|f.
$$
Then we have
$$
|D\sqrt{f}|\geq \sqrt{\frac{|c|}{2}},
$$
at where $f\not=0$. Therefore, we have
$$
\sqrt{f}(s)\geq \sqrt{\frac{|c|}{2}}s
$$
and
$$
f(s)\geq {\frac{|c|}{2}}s^2
$$
 along any minimizing geodesic curve $\gamma(s)$ connecting $x$ and
the fixed point $o$, where $s$ is the arc-length parameter.

Note that using (\ref{B3})we have
$$
|Df|^2(s)=-2cf(x)+R(o)-R(x)\leq -2cf(x)+R(o)\leq Cs^2+R(o).
$$
Hence, for $s>>1$,
\begin{equation}
C_4s\leq |Df|(s)\leq C_5 s. \label{B4}
\end{equation}

\section{$\epsilon$ pinched solitons}

We give a proof of Theorem 3 below. We try to make the proof more
transparent and self-contained.

\emph{Proof of Theorem 3:}
 Recall that
the Ricci curvature of the non-shrinking soliton $(X,g)$ is
$\epsilon$-pinched, i.e., for some $\epsilon>0$ we have that $R>0$
on $X$ and
$$
Rc\geq \epsilon Rg
$$
on $X$. Then using the maximum principle, we know that either $R=0$
on $X$ or $R>0$. If $R=0$ on $X$, then by the pinching condition we
know that $(X,g)$ is Ricci flat.

 Assume that $R>0$ on $X$. Then as
before, the potential function $f$ is locally strictly convex. Since
 $(X,g)$ is a complete non-compact Riemannian manifold, we have that
 $f$ has at most one critical point,i.e., the point where $\nabla
 f=0$. Assume that we have a critical point for $f$, saying that it is $o\in X$.
 Then using (\ref{RR}), we know it is also a critical point of $R$.
 Using (\ref{De}), we know that is the maximum point for $R$. In
 particular, we know that $R$ is a bounded function on $X$,
 saying that $D>0$ is the upper bound.

Using (\ref{RR}) and the $\epsilon$-pinched condition, we have that
$$
-R|\nabla f|^2\leq <\nabla R,\nabla f>=-2Rc(\nabla f,\nabla f)\leq
-\epsilon R|\nabla f|^2.
$$

 Taking a minimizing geodesic curve $\gamma(s)$ connecting $x$ and a
fixed point $o$, where $s$ is the arc-length parameter, we have
\begin{align}
<\nabla f,
\gamma'(s)>|_0^s&=\int_0^sf_{ij}\frac{d\gamma^i}{ds}\frac{d\gamma^j}{ds}ds\nonumber\\
&=\int_0^s(R_{ij}-cg_{ij})\frac{d\gamma^i}{ds}\frac{d\gamma^j}{ds}ds\nonumber\\
&=-cs+\int_0^sR_{ij}\frac{d\gamma^i}{ds}\frac{d\gamma^j}{ds}ds.
\label{PI}
\end{align}
This implies that there is a constant $C_2$ such that
\begin{align*}
<\nabla f, \gamma'(s)> &\geq -cs+\int_0^s\epsilon Rds \geq
-cs+\int_0^1Rds &\geq -cs+C_2\geq C_2
\end{align*}
for $s>>1$.

Using (\ref{B1}) and the pinching condition, we have
 that
 $$
\int_0^sRds\leq C_6.
 $$

Using the pinching condition again, (\ref{PI}) also implies that
$$
<\nabla f, \gamma'(s)>\leq -cs+\int_0^sRds\leq -cs+D.
$$
Therefore, the angle between $\nabla f$ and the gradient of the
distance function from $o$ is almost fixed.

Then, using (\ref{RR}) and the $\epsilon$-pinched condition, we have
for some constant $C_3>0$,
$$
(R^{-1})_s=-R^{-2}<\nabla R, \gamma'(s)> =2R^{-2}Rc(\nabla
f,\gamma'(s)).
$$
Using (\ref{B5}) and (\ref{B4}), we obtain that
\begin{align*}
Rc(\nabla f,\gamma'(s))&=|\nabla f|Rc(\gamma',\gamma')+0(\sqrt{s}) \\
&\geq \epsilon R|\nabla f|+0(\sqrt{s})\geq R(Cs-C),
\end{align*}
we have
$$
(R^{-1})_s\geq 2R^{-1}(Cs-C).
$$
This implies that
$$
(\log R)_s\leq C-Cs
$$
and
$$
R(s)\leq R(o)e^{Cs-Cs^2}.
$$
This implies that $R\to 0$ exponentially as $s\to+\infty$.
\emph{This completes the proof of Theorem 3}.

Theorem 3  tells us that for such $(X,g)$ we have
$$
A=lim sup_{s\to\infty}Rs^2=0.
$$

\end{document}